\documentclass{article}

\usepackage{amssymb,dsfont,latexsym,amsmath}

\newtheorem{thm}{Theorem}

\newtheorem{lem}[thm]{Lemma}
\newtheorem{cor}[thm]{Corollary}

\newcommand\enu[1]{\smallskip\newline\makebox[5mm][l]{\rm(#1)}}

\newcommand\bp{\noindent{\it Proof.}\ }

\begin{document}

\author{Erling St{\o}rmer}

\date{21-9-2010 }

\title{A completely positive map associated with a positive map}

\maketitle
\begin{abstract}
We show that each positive map from $B(K)$ to $B(H)$ is a scalar multiple of a map
of the form $Tr -\psi$ with $\psi$ completely positive.  This is used
to give necessary and sufficient conditions for maps to be $\1C$-positive
for a large class of mapping cones; in particular we apply the results
to k-positive maps. 
\end{abstract}

\section*{Introduction}

In  \cite {SS} we studied several norms on positive maps from $B(K)$
into $B(H)$, where $K$ and $H$ are finite dimensional Hilbert
spaces. These norms were very useful in the study of maps of the form
$Tr-\lambda\psi$, where Tr is the usual trace on $B(K),  \lambda >
0$, and $ \psi$ a completely positive map of $B(K)$ into $B(H)$.  In the
present paper we shall see that every positive map is a positive
scalar multiple of a map of the above form with $ \lambda = 1$, hence
the results in  \cite {SS}  are applicable to all positive maps.  In
particular they yield a simple criterion for some maps to be
k-positive but not k+1-positive.  As an illustration we give a new
proof that the Choi map of $B( \1C^3)$ into itself is atomic, i.e. not
the sum of a 2-positive and a 2-copositive map.

\section*{$\1C$-positive maps}

Let $K$ and $H$ be finite dimensional Hilbert spaces.  We denote by
$ B(B(K),B(H)) $  (resp.$ B(B(K),B(H))^+ $)  the linear (resp. positive
linear) maps of $B(K)$ into $B(H)$. In the case $K=H$ we denote by
$P(H) =   B(B(H),B(H))^+$.  Following \cite {s1} we say a closed cone
$\1C \subset P(H)$ is a \textit{mapping cone} if
$\alpha\circ\phi\circ\beta\in \1C$ for all $\phi\in \1C$ and $\alpha,
\beta\in CP$ - the completely positive maps in $P(H)$.  A map $\phi$
in $ B(B(K),B(H)) $ defines a linear functional $\tilde{\phi}$ on
$B(K)\otimes B(H)$, identified with $B( K\otimes H)$ in the sequel, by
$\tilde{\phi}(a\otimes b) = Tr(\phi(a)b^t)$, where Tr is the usual
trace on $B(H)$ and t denotes the transpose. Let $P(B(K),\1C)$ denote
the closed cone
$$
P(B(K),\1C) = \{a\in B(K\otimes H): \iota\otimes\alpha(a) \geq 0 \ \
\forall\ \  \alpha\in \1C\},
$$
where $\iota$ denotes the identity map on $B(K)$.  Then a map $\phi\in
B((B(K),B(H))$ is said to be $\1C$\textit{-positive} if $\tilde\phi$ is
positive on $P(B(K),\1C)$.  We denote by $\1P_\1C$ the cone of
$\1C$-positive maps.

If $(e_{ij})$ is a complete set of matrix units for $B(K)$ then the
\textit{ Choi matrix} for a map $\phi$ is 
$$
C_\phi = \sum e_{ij}  \otimes \phi(e_{ij} )\in B(K\otimes H).
$$
By \cite{s3} and \cite{s4} the transpose $C_\phi ^t $ of $C_\phi$ is the density
operator for $\tilde\phi$, and by \cite{Ch} $\phi$ is completely
positive if and only if $C_\phi\geq 0$ if and only if $\tilde\phi\geq
0$ as a linear functional on $B(K\otimes H)$.  In the case $\1C = CP$,
$P(CP,B(K)) = B(K\otimes H)^+$, so $\phi$ is CP-positive if and only if
$\phi$ is completely positive. 

If $\1C_1 \subset \1C_2$ are two mapping cones on $B(H)$, then
$P(B(K),\1C_1) \supset P(B(K),\1C_2)$, because if
$\iota\otimes\alpha(a)\geq 0$ for all  $ \alpha\in \1C_2$, then the same inequality
holds for all $\alpha\in \1C_1$. Thus $\tilde\phi\geq 0$ on
$P(B(K),\1C_1) $ implies $\tilde\phi\geq 0$ on $P(B(K),\1C_2), $ so
$\1P_{\1C_1} \subset  \1P_{\1C_2}$. 

Let $\1C$ be a mapping cone on $B(H)$. Let $\1P_C ^o$ denote the
\textit{dual cone of }$\1P_\1C$ defined as 
$$
\1P_C ^o = \{\phi\in B(B(K),B(H)): Tr(C_\phi C_\psi)\geq 0\ \ \forall 
\psi\in \1P_\1C\}.
$$
Thus if $\1C_1 \subset \1C_2$ then $\1P_{\1C_1} ^o\supset \1P_{\1C_2} ^o $.
In the particular case when $\1C\supset CP$ we thus get $\1P_{\1C}^o
\subset  \1P_{CP}^o = CP(K,H)$ - the completely positive maps of
$B(K)$ into $B(H)$.

Following \cite{SS} $\1C$ defines a norm on $B(B(K),B(H))$ by
$$
\parallel\phi\parallel_\1C = sup \{ \mid Tr(C_\phi C_\psi)\mid: \psi\in
\1P_C ^o, Tr(C_\psi)=1\}.
$$
In the special case when  $\1C\supset CP$ it follows from the above
that
$$ 
\parallel\phi\parallel_\1C = sup \mid \rho(C_\phi)\mid,
$$
where the sup is taken over all states $\rho$
on $ B(K\otimes H)$ with density operator $ C_\psi $ with $\psi\in
\1P_C^o$.
Let $\phi\in B(B(K),B(H))$ be a self-adjoint map, i.e. $\phi(a)$ is
self-adjoint for $a$ self-adjoint.  Then $C_\phi$ is a self-adjoint
operator, so is a difference $C_\phi^+ - C_\phi^-$ of two positive
operators with orthogonal supports.  Let $c\geq 0$ be the smallest
positive number such that $c1\geq C_\phi$.  Then $c= \parallel
C_\phi^+\parallel$.  Hence, if $c\neq 0$ there exists a map $\phi_{cp} \in
B(B(K),B(H))$ such that the Choi matrix for $\phi_{cp}$ equals $ 1 -
c^{-1}  C_\phi,$  which is a positive operator. Thus,
if we let  Tr denote the map $x\mapsto Tr(x)1, \phi_{cp}$  is
completely positive, and $c^{-1} \phi = Tr - \phi_{cp}$, since $C_{Tr} =
1$, as is easily shown. Combining the above discussion with \cite{SS},
Prop. 2, we thus have.

\begin{thm}\label{thm1}
Let $\phi$ be a self-adjoint map of $B(K)$ into $B(H)$.  Then if
$-\phi$ is not completely positive,  we have
\enu{i} There exists a completely positive map $\phi_{cp} \in
B(B(K),B(H))$ such that $\parallel C_\phi^+\parallel ^{-1}  \phi = Tr -
\phi_{cp} $.
\enu{ii} If $\1C$ is a mapping cone on $B(H)$ containing CP then
$\phi$ is $\1C$-positive if and only if 
$$
1 \geq \parallel\phi\parallel_\1C = sup\  \rho(C_{\phi_{cp}}),
$$
where the sup is taken over all states $\rho$ on $B(K\otimes H)$  with
density operator $C_\psi$ with $\psi\in\1P_\1C ^o$.
\end{thm}

Note that we did not need to take the absolute value of
$\rho(C_{\phi_{cp}})$ because $C_{\phi_{cp}}\geq 0$ and
$\psi\in\1P_\1C^o \subset CP$. 

We next spell out the theorem for some well known mapping cones.
Recall that a map $\phi$ is \textit{decomposable} if $\phi = \phi_1
+\phi_2$ with $\phi_1$ completely positive and $\phi_2$ copositive,
i.e. $\phi_2 = t\circ \psi$ with $\psi$ completely positive. Also
recall that a state $\rho$ on $B(K\otimes H)$ is a \textit{PPT
  -state} if $\rho\circ(\iota\otimes t)$ is also a state.

\begin{cor}\label{cor2}
Let $\phi\in B(B(K),B(H))$ be a self-adjoint map.  Then we have.
\enu{i} $\phi$ is positive if and only if $\rho(C_{\phi_{cp}})\leq 1$
for all separable states $\rho$ on $B(K\otimes H).$
\enu{ii} $\phi$ is decomposable if and only if
$\rho(C_{\phi_{cp}})\leq 1$ for all PPT-states $\rho$ on $B(K\otimes H).$ 
\enu{iii} $\phi$ is completely positive if and only if
$\rho(C_{\phi_{cp}})\leq 1$ for all states $\rho$ on $B(K\otimes H).$ 
\end{cor}

\bp
(i) That $\phi$ is positive is the same as saying that $\phi$ is
$P(H)-$positive.  Since the dual cone of $P(H)$ is the cone of
separable states (i) follows.

(ii) A state $\rho$ is PPT if and only if its density operator is of
the form $C_\psi$ with $\psi$ a map which is both positive copositive,
see e.g. \cite{s3},Prop.4.  But the dual of those maps is the cone of
decomposable maps, see e.g. \cite{SSZ}.  Thus (ii) follows from the
theorem.

(iii) This follows since the dual cone of the completely positive maps
is the cone of completely positive maps, and that the density operator
for a state is positive, hence the corresponding map $\psi$ is
completely positive.

\section*{k-positive maps}

A map $\phi\in B(B(K),B(H)) $ is said to be \textit{k-positive} if
$\phi\otimes \iota\in B(B(K\otimes L), B(H\otimes L))^+$ whenever $L$
is a k-dimensional Hilbert space.  The k-positive maps in $P(H)$ form
a mapping cone $P_k$ containing $CP$.  Denote by $P_k(K,H)$ the cone of
k-positive maps in $B(B(K),B(H))$.  Then we have ,

\begin{lem}\label{lem3}
With the above notation we have $\1P_{P_k} =P_k(K,H)$.
\end{lem}
\bp
We have $P_{k}^o = SP_k$, the k-superpositive maps in $P(H)$, which
is the mapping cone generated by maps of the form $AdV$ defined by
$AdV(a)= VaV^*$, where $V\in B(H), rank V \leq k$, see
e.g. \cite{SSZ}.  By \cite{s4} the dual cone of $\1P_{P_{k}^o}$ is
given by
$$
\1P_{P_{k}^o}^o = \{ \phi\in B(B(K),B(H)): AdV\circ\phi \in CP(K,H)
\ \forall  V\in B(H), rank V\leq k\}.
$$
By \cite{S},Theorem 3, or \cite{SS},Theorem 2, it follows that
$\1P_{P_{k}^o}^o = P_k(K,H)$. By\cite{s1}, Theorem 3.6, $\1P_{P_k}$ is
generated by maps of the form $\alpha\circ\beta$ with $\alpha\in P_k,
\beta\in CP(K,H)$. Let $AdV\circ\gamma, AdV\in SP_k, \gamma\in
CP(K,H)$ be a generator for $\1P_{P_{k}^o}$.  Then 
$$
Tr(C_{\alpha\circ\beta}  C_{AdV\circ\gamma})= Tr( C_{AdV^*
  \circ\alpha\circ\beta}  C_{\gamma} ) \geq 0,
$$
since $AdV^* \circ\alpha $ is completely positive since $\alpha\in
P_k$ and $rank V \leq k$.  Since the above inequality holds for the
generators of the two cones, it follows that $\1P_{P_k} =
\1P_{P_{k}^o}^o = P_k(K,H)$, completing the proof of the lemma.

\medskip

It follows from the above description of $\1P_{P_k}^o$ that the states with density
operators $C_\psi, \psi\in\1P_{P_k}^o$, are the
same as the vector states generated by vectors in the Schmidt
class $S(k)$, i.e. the vectors $y=\sum_{i=1}^{k} x_i \otimes y_i,
x_i\in K, y_i\in H$, where the $x_i$ and $y_i$ are not necessarily
all $\neq 0$.

\begin{thm}\label{thm4}
Let $\phi\in B(B(K),B(H))^+$.  Then we have.
\enu{i} $\phi$ is k-positive if and only if $sup_{x\in S(k),\parallel
  x\parallel =1} (C_{\phi_{cp}}x,x) \leq 1$.
\enu{ii} Suppose $k< min(dimK, dim H)$, and that there exists a unit
vector $y= \sum_{i=1}^k  x_i\otimes y_i \in S(k)$ 
such that
$y\perp  C_\phi y \notin X\otimes Y$, where $X=span (x_i), Y= span
(y_i)$. 
 Then  $\phi$ is not k+1-positive.
\end{thm}.
In order to prove the theorem we first prove a lemma.

\begin{lem}\label{lem5}
Let $A$ be a self-adjoint operator in $B(K\otimes H)$.  Suppose
$y=\sum_{i=1}^k x_i\otimes y_i$ satisfies $(Ay,y)=1$, and $Ay\notin
X\otimes Y$ with $X,Y$ as in Theorem 4. Then there exist a unit
product vector $x\perp X\otimes Y$ and $s\in (0,1)$ such that
$(A(sx + (1-s^2)^{1/2})y),sx +(1-s^2)^{1/2}y) > 1$.
\end{lem}

\bp
Since $Ay\notin X\otimes Y$ there exists a product vector $x\perp
X\otimes Y$ such that $Re(x,Ay)>0$.  Let $s\in (-1,1)$ and $
t=t(s) = (1-s^2)^{1/2}$, and let $f$ denote the function
$$
f(s) = (A(sx + ty),st + ty) = s^2(Ax,x) + t^2(Ay,y) + 2stRe(Ax,y).
$$
Since $(Ay,y)=1$ we get 
$$
f'(0) = 2(1-s^2)^{1/2}Re(Ax,y) > 0.
$$
Therefore, for $s>0$ and near 0 we have $ (A(sx + ty),st + ty)> f(0) =
1$, proving the lemma.

\medskip

\textit{Proof of Theorem 4}. 

(i) is a direct consequence of Theorem 1, since, as noted in the proof
of Lemma 3, the vector states
$\omega_x$ with $x\in S(k)$ generate the set of states with density
operators $C_\psi$ with $\psi\in \1P_{P_k}^o$.

(ii) By Theorem 1 $C_{\phi_{cp}} = 1
- \parallel C_\phi ^+\parallel^{-1}  C_\phi$, so that $(C_{\phi_{cp}}
y,y) = 1$, using the assumption that $C_{\phi} y\perp y$.  Furthermore 
$C_{\phi_{cp}} y =y -  \parallel C_\phi
^+\parallel^{-1}  C_\phi y.$  Since $ C_{\phi}y\notin X\otimes Y$ 
$C_{\phi_{cp}}y\notin X\otimes Y$.  Thus by Lemma 5 there exist a
unit product vector $x\in X\otimes Y$ and $s,t=(1-s^2)^{1/2} > 0$ such
 that $(C_{\phi_{cp}}(sx +ty),sx + ty) > 1$. Since $sx + ty $ is a
 unit vector in S(k+1), $\phi$
is not k+1-positive by part (i), completing the proof of the theorem.

\medskip
\textit{Example}
We illustrate the above results by an application to the Choi map
$\phi\in B(B(C^3),B(C^3))$ defined by
$$
\phi((x_{ij})) = \begin{bmatrix} 
 x_{11} + x_{33} & -x_{12} & -x_{13} \\
 -x_{21} & x_{11 } + x_{22} & -x_{23}\\   
 -x_{31}  & -x_{32} & x_{22}  +x_{33} 
\end{bmatrix}
$$

\medskip 

We have  $C_{t\circ \phi} = (\iota\otimes t)C_\phi$.  So if $y= x\otimes x$ with
$x = 3^{-1/2} (1,1,1) \in C^3$,  then $(C_\phi y,y) =
(C_{t\circ\phi} y,y) =0$, and $C_\phi y \neq 0 \neq C_{t\circ\phi}
y$.  Hence, by Theorem 4, neither $\phi$ nor $t\circ\phi$ is
2-positive, i.e. $\phi$ is neither 2-positive nor 2-copositive.  Since
$\phi$ is an extremal positive map of $B(C^3)$ into itself by
\cite{CL}, $\phi$ cannot be the sum of a 2-positive and a 2-copositive
map, hence $\phi$ is atomic, a result first proved by Tanahashi and
Tomiyama \cite{TT}, and then extended to more general maps by others,
see \cite{Ha} for references.

$\phi$ can also be shown to be a positive map by a straightforward
argument using Corollary 2.

It should be remarked that the Choi map $\phi$ also yields an example
of a PPT-state on $B(C^3)\otimes B(C^3)$ which is not separable.
Indeed, in \cite{st2} we gave an example of a positive matrix in $A$ in
$B(C^3)\otimes B(C^3)$ such that its partial transpose $t\otimes
\iota(A)$ is also positive, and that $\phi\otimes \iota(A)$ is not
positive.  Then $A$ cannot be of the form $\sum A_i\otimes B_i$ with
$A_i$ and $B_i$ positive, hence the state $\rho(x)=Tr(A)^{-1} Tr(Ax)$
is PPT but not separable. An example of a PPT state on
$B(C^3)\otimes B(C^3)$ which is not separable was later exhibited
by P. Horodecki \cite{1Hor}.

Department of Mathematics, University of Oslo, 0316 Oslo, Norway.

e-mail erlings@math.uio.no


\begin{thebibliography}{999}

 

 \bibitem{Ch}
M.D.Choi, {\em Positive linear maps on complex matrices},
Lin. Alg. Appl. 10 (1975), 285-290.

\bibitem{CL}
M.D.Choi and T.Y.Lam, {\em Extremal positive semifinite forms},
Math.Ann. 231 (1977), 1-18.

\bibitem{Ha}
K-C. Ha, {\em Atomic positive linear maps of matrix algebras},
Publ.RIMS, Kyoto Univ.. 34 (1998), 591-599.

 \bibitem{1Hor}
 P.Horodecki, {\em Separability criterion and inseparable mixed states with positive partial
 transposition}, Physics Letters, A 232(1997), 333-339.

\bibitem{S}
L. Skowronek, {\em Theory of generalized mapping cones in the finite
  dimensional case}, arXiv. 1008.3237(quant-ph)

  \bibitem{SS}
L. Skowronek and E.St{\o}rmer, {\em Choi matrices, norms and
  entanglement associated with positive maps of matrix algebras},
arXiv. 1008.3126(quant-ph)

 \bibitem{SSZ}
L. Skowronek, E.St{\o}rmer and K.Zyczkowski, {\em Cones of positive
  maps and their duality relations}, J.Math.Phys. 50 (2009), 062106.

 \bibitem{s1}
 E.St{\o}rmer, {\em Extension of positive maps into $B(H)$}, J.
 Funct. Anal. 66, No.2 (1986), 235-254.

\bibitem{st2}
E.St{\o}rmer, {\em Decomposable positive maps on C*-algebras},
Proc.Amer.Math.Soc. 86 (1982), 402-404.

 \bibitem{s3}
 E.St{\o}rmer, {\em Separable states and positive maps},
  J.Funct.Anal. 254 (2008),2303-2312.

 \bibitem{s4}
E.St{\o}rmer, {\em Duality of cones of positive maps}, Munster
J.Math. 2 (2009), 299-310.

 \bibitem{TT}
K.Tanahashi and J.Tomiyama, {\em Indecomposable positive maps in
  matrix algebras}, Canad.Math.Bull. 3 (1988),308-317.



 \end{thebibliography}
\end{document}